\documentclass[preprint]{elsarticle}

\usepackage{lineno,hyperref}

\usepackage{amsmath}
\usepackage{amsfonts}
\usepackage[english]{babel}
\usepackage{latexsym}
\usepackage{amssymb,amsthm}
\usepackage{verbatim}
\usepackage{url}
\usepackage{array}

\newtheorem{theorem}{Theorem}[section]

\newtheorem{lemma}[theorem]{Lemma}
\newtheorem{corollary}[theorem]{Corollary}

\newtheorem{proposition}[theorem]{Proposition}
\newtheorem{remark}[theorem]{Remark}

\begin{document}

\begin{frontmatter}

\title{A probabilistic construction of small complete caps in projective spaces}

\author[CorrAuth]{Daniele Bartoli}
\fntext[CorrAuth]{Corresponding author}
\ead{daniele.bartoli@dmi.unipg.it}

\author{Stefano Marcugini}
\ead{gino@dmi.unipg.it}

\author{Fernanda Pambianco}
\ead{fernanda@dmi.unipg.it}

\address{Department of Mathematics and Computer Sciences of Perugia University, Via Vanvitelli 1, Perugia, 06123}




\begin{abstract}
In this work complete caps in $PG(N,q)$ of size $O(q^{\frac{N-1}{2}}\log^{300} q)$ are obtained by probabilistic methods. This gives an upper bound asymptotically very close  to the trivial lower bound $\sqrt{2}q^{\frac{N-1}{2}}$ and it improves the  best known bound in the literature for small complete caps in projective spaces of any dimension. The result obtained in the paper also gives a new upper bound for $l(m,2,q)_4$, that is the minimal length $n$ for which there exists an $[n,n-m, 4]_q2$ covering code with given $m$ and $q$. 
\end{abstract}


\end{frontmatter}


\section{Introduction}

A $n$-cap in an (affine or projective) Galois space  over the finite field with $q$ elements ${\mathbb F}_q$ is a set of $n$ points no three of which are collinear.  A $n$-cap is said to be complete if it is not contained in a $(n+1)$-cap. A plane $n$-cap  is also called a $n$-arc.

A central problem on caps is determining the spectrum of the possible sizes of  complete caps in a given space, see the survey paper \cite{HS2001} and the references therein. As mentioned above, of  particular interest for applications to Coding Theory is the lower part of the spectrum.

For the size of the smallest complete cap in the projective space $PG(N,q)$ of dimension $N$ over ${\mathbb F}_q$, the trivial lower bound is
\begin{equation}\label{trivial_bound}
\sqrt 2 q^\frac{N-1}{2}.
\end{equation}
General constructions of complete caps whose size is close to this lower bound are only known for $q$ even and $N$ odd; see \cite{GDT1991,DGMP2010,Giulietti2007,PS1996}.
When $N$ is even, complete caps of size of the same order of magnitude as $cq^{N/2}$, with $c$ a constant independent of $q$, are known to exist for
both the odd and the even order case, see \cite{DGMP2010,DO2001,Giulietti2007,Giulietti2007b,GP2007}.

If $q$ is odd and $N\equiv 0\pmod 4$, small complete caps
can be obtained via the product method for caps from bicovering plane arcs. It has been shown in \cite{Giulietti2007b} that the cartesian product of a bicovering $k$-arc  $A$ in  the affine plane $AG(2,q)$ and the  cap of size $q^{\frac{N-2}{2}}$ in the affine space $AG(N-2,q)$ arising from the blow-up of a parabola of $AG(2,q^{(N-2)/2})$  is a complete cap in $AG(N,q)$. Via the natural embedding of $AG(N,q)$ in $PG(N,q)$ it is possible to obtain from a complete cap in $AG(N,q)$ a complete cap in $PG(N,q)$ of the same magnitude. 

In \cite{ABGP2013} the authors  obtain  caps of size $(k+1)q^{\frac{N-2}{2}}$ in $AG(N,q)$ when the $k$-arc $A$ is almost bicovering, that is, a complete arc which bicovers all points in $AG(2,q)\setminus A$ but one.

By similar methods, in \cite{ABGP2014} the authors provide new complete caps in $AG(N, q)$ with roughly $q^{(4N-1)/8}$ points, studying both plane cubics with a node and plane cubics with an isolated double point.

In  \cite{ABGP2013} the existence of  complete caps in $AG(N,q)$, $N\equiv 0\pmod 4$, of size of the same order of magnitude as $2pq^{\frac{1}{2}(N-\frac{1}{4})}$, provided that the characteristic $p$ of ${\mathbb F}_q$ is large enough and $\log_pq>8$, is established.

The exact value $t_2(N,q)$ of the minimum size of a complete cap in $PG(N,q)$ is known only for few pairs $(N,q)$: for instance in the case $N=3$, $t_2(3,q)$ is known only for $q\leq 7$; see \cite[Table 3]{DFMP2009}.

In the case $N=3$ according to the survey paper \cite{HS2001}, the smallest known complete caps in $PG(3,q)$, with $q$ arbitrary large, have size approximately $q^{3/2}/2$ and were presented by Pellegrino in 1998 \cite{Pellegrino1998}. However,  Pellegrino's completeness proof appears to present a major gap, and counterexamples can be found; see \cite[Section 2]{BFG2013}.

In this work, the existence of complete caps of size
\begin{equation}\label{obtained_bound}
O(q^{\frac{N-1}{2}} \log ^{300} q)
\end{equation}
in projective spaces $PG(N,q)$, with $N\geq 3$, is established by probabilistic methods. This gives an upper bound asymptotically very close  to the trivial lower bound \eqref{trivial_bound} and it improves the  best known bound in the literature for small complete caps in projective spaces of any dimension. 

In this paper we use  techniques similar to those of  \cite{KV2003} to extend the results from complete arcs in projective plane to complete caps in any projective space of dimension greater than two. In \cite{KV2003} the authors proved the following theorem; see \cite[Theorem 1.1]{KV2003}.

\begin{theorem}
There are positive constants $c$ and $M$ such that the following holds. In every projective plane of order $q\geq M$, there is a complete arc of size at most
$$q^{1/2} log^c q.$$
\end{theorem}

They used new and powerful concentration results (see \cite{KV2000}) and  in particular they proved this stronger result; see \cite[Theorem 1.2]{KV2003}.

\begin{theorem} There are absolute constants $c$ and $M$ such that in any projective plane of order $q\geq M$, one can find an arc with $\Theta(q^{1/2} log^{1/2} q)$ points whose secants cover all but $q^{1/2} log^c q$ points of the plane.
\end{theorem}

Also, they provided an efficient randomized algorithm which produces the desired arc with probability close to $1$ and runs in $\Theta(log^{5/2} q)$ steps, where each step consist of $O(q^4)$ basic operations; see \cite[Section 2.2]{KV2003}.

The main result of this work, in term of complete caps in any projective space $(N,q)$ is the following.

\begin{theorem}\label{MainTheorem}
There are positive constants $c$ and $M$ such that in every projective space of order $q\geq M$ and dimension $N$, there is a complete cap of size at most
$$q^{\frac{N-1}{2}} log^c q.$$
\end{theorem}
Also, as in \cite{KV2003}, a randomized  algorithm to construct the desired complete caps can be easily deduced.

Let $\mathbb{F}_q$ be the finite field with $q$ elements. A $q$-ary linear code $\mathcal{C}$ of length $n$ and dimension $k$ is a $k$-dimensional linear subspace of $\mathbb{F}_q^n$. The Hamming weight $w({\bf x})$ of ${\bf x}$ is the number of nonzero positions in a vector ${\bf x } \in \mathcal{C}$. The minimum distance of $\mathcal{C}$ is defined as
$$d(\mathcal{C}) := \min \{w({\bf x})| {\bf x} \in  \mathcal{C}, {\bf x} \neq 0\}$$
and a $q$-ary linear code of length $n$, dimension $k$, and minimum distance $d$ is denoted as an $[n,k,d]_q$-code. An $[n,k,d]_q$-code can correct at most $t$-errors, where $t=\lfloor \frac{d-1}{2}\rfloor$. The covering radius of $\mathcal{C}$ is the minimum integer $R(\mathcal{C})$ such that for any vector $v \in \mathbb{F}_q^n$ there exists $x\in  \mathcal{C}$ with $w(v-x) \leq R$. An $[n,k,d]_q$-code with covering radius $R$ is denoted by $[n, k, d]_qR$. If $R = t$ then $\mathcal{C}$ is said to be perfect. As there are only finitely many classes of linear perfect codes, of particular interest are those codes $\mathcal{C}$ with $R = t + 1$, called quasi-perfect codes; see  \cite{BLP1998,CHLL1997,CKMS1985}.
The covering density $\mu(\mathcal{C})$, introduced in \cite{CLS1986}, is one of the parameters characterizing the covering quality of an $[n,k,d]_qR$-code $\mathcal{C}$ and it is defined by
$$\mu(\mathcal{C}) =\frac{1}{q^{n-k}}\sum_{i=0}^{R} (q-1)^i \binom{n}{i}.$$
Note that $\mu(\mathcal{C})\geq 1$, and that equality holds when $\mathcal{C}$ is perfect.
Clearly, among codes with the same codimension and covering radius the shortest ones have the best covering density. Therefore  the problem of determining the minimal length $n$ for which there exists an $[n,n-m, d]_qR$-code with given $m, q, d,$ and $R$, has been broadly investigated, see \cite{BPW1989}. Throughout, such minimal length will be denoted as $l(m,R,q)_d$.

In the case  $R = 2$ and $d = 4$, that is, for quasi-perfect linear codes that are both $1$-error correcting and $2$-error detecting, the columns of a parity check matrix of an $[n, n- m, 4]_q2$-code can be considered as points of a complete $n$-cap in the finite projective space $PG(m- 1, q)$. For this reason these codes have been investigated also through their  connection with Projective Geometry;  see e.g. \cite{GP2007}.

The construction of  complete caps presented in this paper gives, via the connection between Coding Theory and Projective Geometry, the following upper bound on $l(N+1,2,q)_4$.
\begin{theorem}\label{MainTheoremCodici}
There exists a positive constants  $M$ such that for every  $q\geq M$
$$l(N+1,2,q)_4=O(q^{\frac{N-1}{2}} \log ^{300} q).$$
\end{theorem}

The paper is organized as follows. In Section \ref{Algorithm} the algorithm based on the nibble method for the construction of the compete caps is described. In Section \ref{Notations} the notations used in this work are presented. In Section \ref{SMainLemma} the main features of the two phases of the algorithm are listed and the Main Theorem is proven via these properties. In Section \ref{Instruments} all the principal probabilistic instruments are collected. In Sections \ref{Phase1} and \ref{Phase2} all the properties used in the proof of the main Theorem are proven.

\section{Algorithm}\label{Algorithm}
\subsection{The nibble method}
The algorithm used to construct complete caps is the same as in \cite{KV2003}.  This algorithm is based on the nibble method and in the following we recall the main features of this method. For a more detailed introduction see \cite[Section 2.1]{KV2003}.

Usually to construct a cap one can proceed selecting the points one by one (for instance see \cite{BDFMP2013,BDFMP2014b}) adding at each step a point not lying on bisecants of the cap at the present step. The choice can be done in a complete random way, taking the first point in a certain order (see \cite{BDFMP2014b}), or selecting the point which maximizes an objective function (see \cite{BDFMP2013}). Usually the results obtained in this way are good if we are interested in small complete caps; see \cite{BDFMP2013,BDFMP2014b,BDFKMP2013b,BDFMP2012,BDFMP2014} and the reference therein.

The nibble method, or dynamic random construction using nibble (DRC), is an approximated version of the previous algorithm. Instead of selecting one element at each step, DRC randomly and independently chooses points, not selected yet, with a certain probability so that a bunch of elements are chosen together. This is called a ``nibble'' and its size is the number of chosen elements or  its expectation.
The first problem, when selecting a set of points instead of a single point at each step, is that the new set obtained by the union of the nibble and the cap at the previous step could be not a cap. In order to avoid this, only a subset of the nibble is added to the cap. Usually, each point in the nibble (we will refer to the points in the nibble as \emph{chosen} points) which does not lie on bisecants of the cap or of the nibble itself is \emph{selected}, i.e. is chosen to be added to the constructing cap.
As not all the chosen points are selected, some points could  unnecessarily be eliminated, but in this way the set of remaining elements available for the subsequent steps can be well-understood.
The size of a nibble cannot be too big, since in this case  many of the chosen points contribute to at least one conflict (that is lies on some bisecants), and it would be hard to predict the structure of the cap at each step and/or too many elements would unnecessarily be discarded.
Therefore the size of the nibble should be small enough to have that most of the chosen elements do not cause any conflict and only few elements would be unnecessarily discarded.
For instance, if we choose $\theta q^{\frac{N-1}{2}}$ random points from a projective space of order $q$, each chosen point causes a conflict  with probability at most
$$(q^{N-1}+q^{N-2}+\ldots+q+1)\binom{q}{2} \left(\theta q^{-\frac{N+1}{2}}\right)^2 \leq \theta^2.$$
If $\theta$ is small enough, then most of the chosen points do not cause any conflict.

\subsection{The algorithm}
The cap is constructed in the following way. At the beginning, the starting cap $A_0$ is empty. Let $\Omega_0=S_0$ be the set of all the points of the projective space $PG(N,q)$. Roughly speaking, at each step $\Omega_i$ is essentially the set of points which are not covered by the cap  $A_i$, while $S_i$ is a subset of $\Omega_i$. At the each step $i$, a random subset $B_i \subset S_i$ is selected, choosing each point from $S_i$ independently with the same probability $p_i$. The set  $B_i$ is the \emph{nibble} and only a subset of $B_i$ is added to  $A_i$ to obtain the new cap $A_{i+1}$. This subset $M_i \subset B_i$ is the set of the points not causing any conflict.
At the subsequent step,  $\Omega_{i+1}$ is obtained from $\Omega_i$ by deleting all the points covered by the secants of $A_{i+1}$ or in $B_i$, while $S_{i+1}$ is obtained by deleting from $S_i$ the points covered by the secants of $A_{i+1}$ or in $B_i$ plus a few more points, chosen randomly: in this way certain structural properties of the $S_i$'s are preserved.
The process is repeated until all but $q^{\frac{N-1}{2}} \log^c q$ points are covered by the secants of the current arc. In the following we set $\theta=\log^{-2} q$.

The algorithm acts as follows. At each step we use three different ``sub-phases'': {\bf choose}, {\bf delete}, and {\bf compensate}.

{\bf Start} : $\Omega_0$, $S_0$ are both $PG(N,q)$, $A_0=\emptyset$. We also consider the  quantities $a_i$ and $b_i$. At the beginning $a_0=0$ and $b_0=1$, while at the $i$-th step $$a_i=\frac{|A_i|}{q^{\frac{N-1}{2}}} \qquad \textrm{ and } \qquad b_i=\frac{|S_i|}{q^{N}+q^{N-1}+\cdots+q+1}.$$

{\bf Choose} : At each step a point $v$ in $S_i$ is chosen with probability $p_i=\theta(b_i q^{\frac{N+1}{2}})^{-1}$. The set of all the chosen points is $B_i$. A point $x$ in $B_i$ is \emph{good} in $A_i \cup B_i$ if there are no two points in $A_i \cup B_i$ collinear with $x$. The set $M_i$ is the set of all the good points. So $A_{i+1}=A_i \cup M_i$.

{\bf Delete} : Delete from $\Omega_{i}$ all the points in bisecants of $A_{i+1}$ or in $B_i$. Let $D_{i}$ be the set of deleted points and, if $v \in \Omega_i$, let $P_{i}(v) =Pr(v \in D_i)$. Let $P_i^{u}$ and $P_i^{\ell}$ be the upper and the lower bounds for these probabilities.

{\bf Compensate} : $S_{i+1}$ is obtained from $S_i$ deleting the points of $D_i$ and independently the points of $S_i$ with probability
$$P_i^{com}(v)=\frac{P_i^u-P_i(v)}{1-P_i(v)}.$$
Let $R_i$ be the set of the removed points, then

$$\Omega_{i+1}=\Omega_{i}\setminus D_i, \quad S_{i+1}=S_i \setminus (D_i \cup R_i), \quad A_{i+1}=A_i \cup M_i$$
and
$$ b_{i+1}=b_{i}(1-P_i^{u})=\prod_{j=1}^{i}(1-P_{j}^{u}), \qquad b^{\prime}_{i+1}=b^{\prime}_{i}(1-P_i^{\ell})=\prod_{j=1}^{i}(1-P_{j}^{\ell}).$$

{\bf Stop} : The algorithm stops after $K$ steps if $K$ is the first integer such that
$$b_{K}\leq q^{-\frac{N+1}{2}}\log ^{c} q,$$
for some constant $c$ (we will set $c=300$, as in \cite{KV2003}).

\begin{remark}\label{Remark_i}
All the properties shown in the following are proven for all $i\leq o(\log^{3} q)$.
\end{remark}

The importance of the sub-phase of Compensation is explained in the following remark.
\begin{remark}\label{RemarkComp}
The operation of ``compensation'' is made in order to give the same probability to the points in $S_i$ to be in $S_{i+1}$. In fact, if $p=P_i(v)$, then
$$P(v \notin S_{i+1}| v \in S_i)=p+(1-p)\frac{P_i^u-p}{1-p}=P_i^u.$$
So,
$$\mathbb{E}(|S_{i+1}|)=|S_i|(1-P_i^u).$$
\end{remark}

In order to prove the properties on the nibble and on the cap at every step, we divide the search into two phases.
The first phase continues when
$$b_i\geq \frac{\log^{c_1} q}{q},$$
while the second phase continues when
$$ \frac{\log^{c} q}{q^{\frac{N+1}{2}}}\leq b_i \leq \frac{\log^{c_1} q}{q}.$$

\section{Notations}\label{Notations}

Let $\ell$ be a line and $u,v$ points of $PG(N,q)$. Also, if $a,b,c$ are points of the space, then $[abc]$ indicates that the three points are collinear, while $(ab)$ is the line through $a$ and $b$. In the following we will use the following sets.
\begin{itemize}
\item $A_i(v)=\{x \in S_i \setminus \{v\} | \exists u \in A_i : [xuv]\}$ is the set of all the points of $S_i$ belonging to the cone of lines through $v\in \Omega_i$ and the points of $A_i$; $A_i^{\prime}(v)=\{x \in \Omega_i \setminus \{v\} | \exists u \in A_i : [xuv]\}$ is the set of all the points of $\Omega_i$ belonging  to the cone of lines through $v\in \Omega_i$ and the points of $A_i$; $A_i(u_1,\ldots,u_m)=\bigcap_{j=1}^{m}A_i(u_j)$;
\item $B_i(v)=\{x \in S_i \setminus \{v\} | \exists u \in B_i : [xuv]\}$ is the set of all the points of $S_i$ belonging to the cone of lines through $v\in \Omega_i$ and the points of $B_i$; $B_i(u_1,\ldots,u_m)=\bigcap _{j=1}^{m}B_i(u_j)$;
\item $S_i(\ell)= S_i \cap \ell$: note that if $\ell$ is a bisecant of the cap then $S_i(\ell)$ is empty; $S_i(\ell,u)=S_i(\ell)\cap A_i(u)$; $S_i(\ell,u,v)=S_i(\ell)\cap A_i(u)\cap A_i(v)$;
\item $\Omega_i(\ell)= \Omega_i \cap \ell$; $\Omega_i(\ell,u)=\Omega_i(\ell)\cap A_i^{\prime}(u)$; $\Omega_i(\ell,u,v)=\Omega_i(\ell)\cap A_i^{\prime}(u) \cap A_i^{\prime}(v)$;
\item $T_i(v)=\{\{x,y\} | x,y \in S_i\setminus \{v\} , [vxy]\}$, the set of the unordered pairs of points of $S_i$ collinear with $v$; $T_i^{\prime}(v)=\{(x,y) | x,y \in S_i\setminus \{v\} , [vxy]\}$, the set of the ordered pairs of points of $S_i$ collinear with $v$.
\end{itemize}
From the above definitions the following hold:
$$|T_{i}(v)|=\sum_{\ell \ni v }\binom{|S_i(\ell)|-1}{2}; \qquad |S_i|=\sum_{\ell} \frac{|S_i(\ell)|}{q^{N-1}+\cdots+q+1};$$
$$|\Omega_i|=\sum_{\ell} \frac{|\Omega_i(\ell)|}{q^{N-1}+\cdots+q+1}; \qquad |A_{i}(v)|=\sum_{{\scriptsize \begin{array}{c}\ell \ni a,v \\ a \in A_{i}\\ \end{array}} }(|S_i(\ell)|-1).$$
Also, given two functions $f(q), g(q)$, as usual, we say that $g(q)=o(f(q))$, $g(q)=O(f(q))$, $g(q)=\Theta(f(q))$, if
$$\lim_{q \to +\infty} \frac{g(q)}{f(q)}=0,\quad \lim_{q \to +\infty} \frac{g(q)}{f(q)}\leq K \in \mathbb{N}, \quad \lim_{q \to +\infty} \frac{g(q)}{f(q)}=K\in \mathbb{N},$$
respectively.
In the following we say that an event $A$ ``occurs with very high probability'' if $Pr(\overline A)\leq e^{-w\log q}$ for some constant $w \in \mathbb{N}$. We consider in all the proofs that the dimension $N$ of the space is \emph{greater than} $2$.

\section{Main Lemma}

In the following lemma all the properties which hold at each step are summarized.

\begin{lemma}\label{Main Lemma}\label{SMainLemma}
The following holds for any $i\leq o(\log^3 q)$.
\begin{itemize}
\item {\bf Phase 1}
\begin{itemize}
\item {\bf Primary Properties}
\begin{eqnarray}
\theta q^{\frac{N-1}{2}}(1-o(1))\leq |M_i| \leq \theta q^{\frac{N-1}{2}}(1+o(1)) \qquad |B_i| \leq 2 \theta q^{\frac{N-1}{2}}\label{1}\\
b_{i+1}q(1-(i+1)\log^{-13}q)\leq |S_{i+1}(\ell)|\leq b_{i+1}q(1+(i+1)\log^{-13}q)\label{2}\\
|\Omega_{i+1}(\ell)|\leq b_i^{\prime}q(1+(i+1)\log^{-13}q)\label{3}
\end{eqnarray}
\item {\bf Secondary Properties}
\begin{eqnarray}
|S_{i+1}(\ell,v)|\leq 8(i+1)a_{i+1}b_{i+1}\sqrt{q}+(i+1)\log^{40}q\label{4}\\
|S_{i+1}(\ell,u,v)|\leq (i+1)\log^{4}q\label{5}\\
|A_{i+1}(u,v)|\leq (i+1)b_{i+1}q+(i+1)\log^{40}q\label{6}\\
|\Omega_{i+1}(\ell,v)|\leq 8(i+1)a_{i+1}b^{\prime}_{i+1}\sqrt{q}+(i+1)\log^{40}q\label{9}\\
|\Omega_{i+1}(\ell,u,v)|\leq (i+1)\log^{4}q\label{10}
\end{eqnarray}
\end{itemize}

\item {\bf Phase 2}
\begin{eqnarray}
\theta q^{\frac{N-1}{2}}(1-o(1))\leq |M_i| \leq \theta q^{\frac{N-1}{2}}(1+o(1)) \qquad |B_i| \leq 2 \theta q^{\frac{N-1}{2}}\label{1bis}\\
\frac{1}{2}b_{i+1}^2q^{N+1}(1-3(i+1)\log^{-13}q)\leq |T_{i+1}(v)|\leq \frac{1}{2} b_{i+1}^2q^{N+1}(1+3(i+1)\log^{-13}q)\label{2bis}\\
a_{i+1}b_{i+1}q^{\frac{N+1}{2}}(1-O(i\theta^2) )\leq |A_{i+1}(v)|\leq a_{i+1}b_{i+1}q^{\frac{N+1}{2}}(1+O(i\theta^2) ) \label{3bis}\\
(1-\log^{-10}q)^{i+1}b_{i+1}q^{N}\leq |S_{i+1}|\leq (1+\log^{-10}q)^{i+1}b_{i+1}q^{N}\label{4bis}\\
|\Omega_{i+1}|\leq (1+\log ^{-10} q)^{i+1} q^{N} b^{\prime}_{i+1}\label{5bis}
\end{eqnarray}
\end{itemize}
\end{lemma}

\begin{remark}\label{Additional Properties}
From the properties of the first phase, recalling that $p_i=\frac{\theta}{b_iq^{\frac{N+1}{2}}}$, we have

$$|A_i(v)|=|A_i||S_i(\ell)|, \qquad |T_i(v)|\simeq \frac{1}{2}(q^{N-1}+q^{N-2}+\cdots+q+1)|S_i(\ell)|^2$$
and therefore
\begin{eqnarray}
\frac{1}{2}b_{i+1}^2q^{N+1}(1-3(i+1)\log ^{-13} q)\leq |T_{i+1}(v)|\leq \frac{1}{2}b_{i+1}^2q^{N+1}(1+3(i+1)\log ^{-13} q)\label{11}\\
a_{i+1}b_{i+1}q^{\frac{N+1}{2}}(1-2(i+1)\log^{-13} q )\leq |A_{i+1}(v)|\leq a_{i+1}b_{i+1}q^{\frac{N+1}{2}}(1+2(i+1)\log^{-13} q ) \label{12}\\
(1-\log^{-10}q)^{i+1}b_{i+1}q^{N}\leq |S_{i+1}|\leq (1+\log^{-10}q)^{i+1}b_{i+1}q^{N}\label{13}\\
|\Omega_{i+1}|\leq (1+\log^{-10} q)^{i+1} q^{N} b^{\prime}_{i+1}\label{14}
\end{eqnarray}
\end{remark}

\section{Proof of the Main Theorem via the Main Lemma}
In this section we prove the Main Theorem using the properties listed in the Main Lemma. The proofs of Lemmas \ref{Lemma 3.2} and  \ref{Lemma 3.3} are the same as  \cite[Lemma 3.2]{KV2003} and \cite[Lemma 3.3]{KV2003} and therefore omitted.

\begin{lemma}\label{Lemma 3.2}
We have
$$P_i^u\leq p_i+p_i\max_{v} |A_i(v)|+p_i^2\max_{v}|T_i(v)|,$$
where the maximum is taken over the set $\Omega_i$.
\end{lemma}
The following estimates hold:
\begin{eqnarray*}
p_i=\frac{\theta}{b_i q^{\frac{N+1}{2}}}\leq \frac{\theta}{\log^c q}=log^{-302} q\leq \log^{-4}q.\\
p_i^2\max_{v}|T_i(v)|\simeq  \frac{\theta^2}{b_i^2q^{N+1}}b_i^2 q^{N+1}=\theta^2=\log^{-4}q.\\
\log^{-4}q <p_i\max_{v}|A_i(v)|\simeq \frac{\theta}{b_iq^{\frac{N+1}{2}}}ib_i\theta q^{\frac{N+1}{2}}=i\theta^2\leq \log^{-1}q.
\end{eqnarray*}

\begin{lemma}\label{Lemma 3.3}
We have
$$P_i^l\geq p_i\min_{v}|A_i(v)|-2p_i^2\max_{v} |A_i(v)|^2-p_i^3\max_{v}|A_i(v)|\max_{v}|T_i(v)|,$$
where the maximum is taken over the set $\Omega_i$.
\end{lemma}

The following estimates hold:
\begin{eqnarray*}
p_i\min_{v}|A_i(v)|\simeq \frac{\theta}{b_iq^{\frac{N+1}{2}}}ib_i\theta q^{\frac{N+1}{2}}=i\theta^2.\\
2p_i^2\max_{v} |A_i(v)|^2\simeq 2 i^2\theta^4=o( i \theta^2).\\
p_i^3\max_{v}|A_i(v)|\max_{v}|T_i(v)|\simeq  i\theta^2 \frac{\theta^2}{b_i^2q^{N+1}}b_i^2 q^{N+1}=i\theta^4=o(i\theta^2).
\end{eqnarray*}

Let $K$ be the first integer such that
$$b_K< \frac{\log^c q}{q^{\frac{N+1}{2}}},$$
that is $K$ is the total number of steps of the algorithm.

\begin{lemma}\label{Lemma 3.5}
Suppose that Lemma \ref{Main Lemma} holds. Then
 $$K=\Theta(\theta^{-1}\log^{1/2}q)=\Theta(\log^{5/2} q).$$
\end{lemma}
\proof
By the definition of $K$
$$\prod_{i=0}^{K-1}(1-P_i^u)=b_{K}\leq q^{-\frac{N+1}{2}}\log ^c q\leq b_{K-1}=\prod_{i=0}^{K-2}(1-P_i^u).$$
Suppose $K\geq L=\left \lceil\sqrt{6N}\log^{5/2}q\right\rceil$, then
\begin{equation}\label{b_L}
b_{L}\geq b_{K-1}\geq q^{-\frac{N+1}{2}}\log ^c q
\end{equation}
and
$$b_{L} =\prod _{i=1}^{L-1}(1-P_i^u)\leq \prod _{i=1}^{L-1}(1-P_i^{\ell}).$$
From \eqref{b_L} we have
$$-\frac{N+1}{2}\log q+c\log \log q =\left(-\frac{N+1}{2}+\frac{\log \log q}{\log q}\right)\log q=$$
$$=\left(-\frac{N+1}{2}+o(1)\right)\log q  \leq \sum_{i=1}^{L-1}\log (1-P_i^{\ell}).$$
From Lemma \ref{Lemma 3.3} we have
$$P_i^{\ell}\geq \frac{1}{2}i \theta ^2,$$
so
$$\log (1-P_i^{\ell})\leq \log \left(1-\frac{1}{2}i\theta^2\right)\leq- \frac{1}{3}i \theta^2.$$
Then
$$\sum_{i=1}^{L-1}\log (1-P_i^{\ell})\leq -\frac{1}{3}\sum_{i=1}^{L-1} i \theta^2\leq -\frac{1}{7}L^2\theta^2 \leq -\frac{6}{7}N\log q$$
which is impossible, then $K< L=\left\lfloor \sqrt{6N}\log^{5/2}q\right\rfloor$.
Suppose that  $K\leq L=\left\lfloor \sqrt{\frac{N}{6}}\log^{5/2}q\right\rfloor$, then
\begin{equation}\label{b_L_seconda}
\prod _{i=1}^{L}(1-P_i^u)=b_{L+1}\leq b_{K}\leq q^{-\frac{N+1}{2}}\log ^c q
\end{equation}
and then
$$\sum_{i=1}^{L}\log (1-P_i^u)\leq \left(-\frac{K+1}{2}+o(1)\right)\log q$$
From Lemma \ref{Lemma 3.2} we have  $P_i^{u}\leq 2i\theta^2.$
Then
$$\log (1-P_i^{\ell})\geq \log \left(1-2i\theta^2\right)\geq- 3i \theta^2.$$

$$-\frac{K}{2}\log q  \leq \sum_{i=1}^{L}\log (1-P_i^u)\leq  \left(-\frac{K+1}{2}+o(1)\right)\log q,$$
which is impossible, therefore $K>\sqrt{\frac{N}{6}}\log^{5/2}q$ and $K=\Theta(\log^{5/2} q)$.
\endproof

\begin{corollary}\label{corA_N}
Suppose that Lemma \ref{Main Lemma} holds, then
$$|A_K| =\Theta (q^{\frac{N-1}{2}}\log^{1/2} q).$$
\end{corollary}
\proof
By Lemma \ref{Lemma 3.5} the number of steps is $\Theta(\log^{5/2} q)$; since $|A_K|=\sum_{i=1}^{K}|M_i|$, then
$$|A_K|=\theta q^{\frac{N-1}{2}}\Theta(\log^{5/2} q)=\Theta (q^{\frac{N-1}{2}}\log^{1/2} q).$$
\endproof

\begin{lemma}\label{Lemma 3.8}
Suppose that Lemma \ref{Main Lemma} holds, then for every $i$
$$P_i^u -P_i^{\ell}=O(\log^{-3}q).$$
\end{lemma}
\proof
By Lemmas \ref{Lemma 3.2} and \ref{Lemma 3.3} we have
$$P_i^u-P_i^{\ell}\leq p_i(1+\max_v |A_i(v)|-\min_v |A_i(v)|)+p_i^2\max_v|T_i(v)|+$$
$$+2p_i^2\max_v |A_i(v)|^2+p_i^3\max_v |A_i(v)|\max|T_i(v)|.$$
The following estimations hold:
$$p_i^3\max_{v}|A_i(v)|\max_{v}|T_i(v)|=O(i\theta^4)=O(\log ^{-11/2} q);$$
$$2p_i^2\max_v |A_i(v)|^2=O(2 i^2\theta^4)=O(\log^{-3}q);$$
$$p_i^2\max_v|T_i(v)|=O(\theta^2)=O(\log^{-4} q).$$
By Property \eqref{3bis}
$$\max_v |A_i(v)|-\min_v |A_i(v)|=O(i\theta^2a_{i+1}b_{i+1}q^{\frac{N+1}{2}}),$$
then
$$p_i(\max_v |A_i(v)|-\min_v |A_i(v)|)=O(i^2\theta^4)=O(\theta^2\log q)=O(\log^{-3}q),$$
and finally
$$p_i=O(\log^{-2-c}q).$$
\endproof

Using the previous lemmas we are able to prove the main theorem.
\begin{theorem}[Main Theorem]
There are positive constants $c$ and $M$ such that in every projective space of order $q\geq M$ and dimension $n$, there is a complete cap of size at most
$$q^{\frac{N-1}{2}} log^c q.$$
\end{theorem}
\proof Suppose that Lemma \ref{Main Lemma} holds.
By Corollary \ref{corA_N}, $|A_{N}|=\Theta (q^{\frac{N-1}{2}}\log^{1/2}q)$.
By \eqref{14}
$$|\Omega_{K}|\leq (1+\log^{-10} q)^{K} q^{N} b^{\prime}_{K}=(1+o(1))q^{N}\prod_{i=0}^{K-1}(1-P_i^{\ell})$$
and $$|S_{K}|=\Theta(b_{K}q^{N}).$$
Then
$$\frac{|\Omega_{K}|}{|S_{K}|}=O\left(\prod_{i=0}^{K-1}\frac{1-P_i^{\ell}}{1-P_i^{u}} \right).$$
By Lemma \ref{Lemma 3.8} we have that
$$\log \left(\frac{1-P_i^{\ell}}{1-P_i^{u}}\right)=O(P_i^u-P_i^\ell)=O(\log^{-3}q).$$
Then
$$\log \frac{|\Omega_{K}|}{|S_{K}|} =O\left(\sum_{i=0}^{K-1} \log \left(\frac{1-P_i^{\ell}}{1-P_i^{u}}\right) \right)=O(K\log^{-3}q)=O(\log^{-1/2} q)=o(1),$$
and 
$$|\Omega_{K}|=e^{o(1)}|S_{K}|.$$
By Property \eqref{4bis}
$$|S_K|=O(b_Kq^{N})=O(q^Nq^{-\frac{N+1}{2}}\log ^c q)=O(q^{\frac{N-1}{2}}\log ^c q),$$
then
$$|\Omega_{K}|=O(q^{\frac{N-1}{2}}\log ^c q).$$
Also,
from Properties \eqref{1} and \eqref{1bis},
$$|B_i| =O(\theta q^{\frac{N-1}{2}}).$$
Since after the $K$-th step the uncovered points are a subset of $\Omega_K \cup \left( \cup_{i=1}^{K}B_i\setminus M_i\right)$, the size of this  set is $O(q^{\frac{N-1}{2}}\log ^c q)$ and therefore  there exists a complete cap of size
$$O(q^{\frac{N-1}{2}}\log ^c q).$$
\endproof

\section{Instruments}\label{Instruments}
In this section we summarize all the instruments from \cite{KV2003} which will be used for the proofs of the properties of the first and the second phase in the Main Lemma.

One of the most important tool used in this paper is the polynomial method; see \cite[Section 4.2]{KV2003}. In general, when using the probabilistic method, a key point is to show that a particular random variable is strongly concentrated around its mean.
In \cite{KV2003} the authors pointed out that none of the existing (classical) tools was sufficiently strong to prove the fundamental  properties in the nibble process; therefore they presented two concentration results (Lemmas $4.1$ and $4.2$ in \cite{KV2003}) we also use in this paper.

Let $t_1,\ldots,t_r$ be $r$ independent binary random variables with expectation $p_i$. They indicate if the $i$-th point is chosen or not. Let $Y$ be a function depending on $t_1,\ldots,t_r$; then for $i=1,\ldots,r$ let 
$$C_i(v) =\left| \mathbb{E}\left(Y(v^{(1)})-Y(v^{(0)})|t_1,\ldots,t_{i-1}\right)\right|,$$
where $$Y(v^{(1)})=Y(t_1,\ldots,t_{i-1},1,t_{i+1},\ldots,t_{r})$$ and $$Y(v^{(0)})=Y(t_1,\ldots,t_{i-1},0,t_{i+1},\ldots,t_{r}).$$ $C_i(v)$ is called the conditional average effect of the random variable $t_i$ when $t_1,\ldots,t_{i-1}$ are given.
Let $L\subset S_i$ be a subset of size $|L|\geq \log^{100}q$ and $L^{\prime}=L \cap S_{i+1}$. Given a point $j$, define
$$A(L,j)=\{t \in L | (tj)\cap A_i\neq \emptyset\} \; \textrm{ and } \; a(L)=\max \left\{ \max_{j \in \Omega_i}|A(L,j)|,|L|\log^{-100} q\right\}.$$
In the following we list some of the results of \cite{KV2003} we will use in Sections \ref{Phase1} and \ref{Phase2}. The proofs are all very similar to the ones in \cite{KV2003}. Also, in some cases, they are precisely the same (see e.g. Lemmas \ref{Lemma 4.7} and \ref{Lemma 4.2} and Corollary \ref{Corollary 4.8}) and therefore they are omitted.
\begin{lemma}[Lemma 4.2 in \cite{KV2003}]\label{Lemma 4.2}
Let $H$ be a hypergraph with vertex $\mathcal{V}(H)=\{1,2,\ldots,r\}$ and edges $\mathcal{E}(H)$ such that each edge has at most $k$ vertices. Let $w(e)$ the weight of a single edge $e$. Suppose $t_1,\ldots, t_r$ are independent random variables (binary random variables with expected value $p_i$). Consider
$$Y_{H}=\sum_{e \in \mathcal{E}(H)} w(e) \prod _{s \in e} t_s.$$
Let
$$\mathbb{E}_{i}(Y)=\max_{A\subset \mathcal{V}(H) | |A|=i}\mathbb{E}(Y_{H_A}),$$
where $Y_{H_A}$ is the $A$-truncated subgraph of $H$, that is the hypergraph with vertices $\mathcal{V}(H)\setminus A$,  $\mathcal{E}(H_A)=\{B \subset \mathcal{V}(H_A), B \cup A \in \mathcal{E}(H)\}$ and $w(B)=w(B\cup A)$.
Then there exist positive numbers $c_k$, $d_k$ depending only on $k$ so that for any  $\lambda>0$
$$Pr\left(|Y-\mathbb{E}(Y)|\geq (E E^{\prime})^{1/2} \lambda^{k})\right)\leq d_k e^{-\lambda+k\log n}.$$
\end{lemma}

Let $L\subset S_i$ (or $L\subset \Omega_i$). For any point $j$, let
$$A(L,j)=\{t \in L | (tj)\cap A \neq \emptyset\}.$$

\begin{lemma}[Lemma 4.4 in \cite{KV2003}]\label{Lemma 4.4}
With very high probability,
$$C_k(v)=o(\log q) a(L).$$
\end{lemma}
\proof
The proof is the same as in \cite[Lemma 4.4]{KV2003}, since $p_i|A_i(j)|\leq 2a_i\theta=o(1)$ and $2b_i q p_i\leq q^{-\frac{N-1}{2}}\leq log^{-102} q$, for $q$ large enough.
\endproof

\begin{lemma}[Lemma 4.5 in \cite{KV2003}]\label{Lemma 4.5}
With very high probability,
$$\sum_{k=1}^{|S_i|} p_i C_k(v)=o(\log q) |L|.$$
\end{lemma}
\proof
The proof is the same as in \cite[Lemma 4.5]{KV2003}, since $p_i\max_j |A_i(j)|=o(1)$ and $2b_i q p_i\leq q^{-\frac{N-1}{2}}\leq log^{-102} q$, for $q$ large enough.
\endproof

\begin{lemma}[Lemma 4.7 in \cite{KV2003}]\label{Lemma 4.7}
Let $L\subset \Omega$. With very high probability,
$$|L^{\prime}-\mathbb{E}(L^{\prime})|\leq a (L)^{1/2}|L|^{1/2}log^{5} q.$$
\end{lemma}

\begin{corollary}[Corollary 4.8 in \cite{KV2003}]\label{Corollary 4.8}
Let $K$ be a fixed positive constant. For any set $L$ such that $a(L)\leq \frac{|L|}{\log^{2(K+5)}q}$, with very high probability,
$$|L^{\prime}-\mathbb{E}(L^{\prime})|\leq |L|\log^{-K} q.$$
\end{corollary}

The following lemma states that, given two points $x,y \in S_i$, even if the events $x \in S_{i+1}$ and $y \in S_{i+1}$ are not independent, they are ``almost independent''.
\begin{lemma}[Lemma 5.3 in \cite{KV2003}]\label{Lemma 5.3}
Let $x,y \in S_i$, then
$$|Pr(x,y\in S_{i+1})-Pr(x\in S_{i+1})Pr(y \in S_{i+1})|=o(\log^{-13}q).$$
\end{lemma}

\section{First Phase}\label{Phase1}

\begin{proposition}\label{Property1}
Property \eqref{1} holds.
\end{proposition}
\proof
By induction on $i$.
\begin{itemize}
\item $|B_i|\leq 2\theta q^{\frac{N-1}{2}}$ : every point is chosen with the same probability $p_i=\frac{\theta}{b_iq^{\frac{N+1}{2}}}$ from the set $S_i$, then, since by induction $|S_i|\simeq(1+\log^{-10} q)^i b_i q^N$ (see Property \eqref{13}), with very high probability
$$|B_i|\simeq p_i |S_i|=\frac{\theta}{b_iq^{\frac{N+1}{2}}} (1+\log^{-10} q)^i b_i q^N=(1+o(1))\theta q^{\frac{N-1}{2}}\leq 2\theta q^{\frac{N-1}{2}}.$$
\item $\theta q^{\frac{N-1}{2}}(1-o(1))\leq |M_i| \leq \theta q^{\frac{N-1}{2}}(1+o(1))$: to show this, we get a bound on $U_i=B_i \setminus M_i$. Using the polynomial method (see \cite[Subsection 4.2]{KV2003}), consider the variables
$$\{t_j| j \in S_i\},$$
such that $t_j=1$ if $j\in B_i$. If a point is in $B_i$ but not in $M_i$, then  there exists a pair of points $j^{\prime}, j^{\prime\prime}$ collinear with $j$ such that either $j^{\prime},j^{\prime\prime} \in B_i$ or $j^{\prime} \in B_i$ and $j^{\prime\prime} \in A_i$. In the first case
$$\sum_{j^{\prime},j^{\prime\prime}| [j,j^{\prime},j^{\prime\prime}]}t_{j^{\prime}}t_{j^{\prime\prime}}> 0,$$
while in the second case
$$\sum_{j^{\prime}\in A_i(j)}t_{j^{\prime}}> 0.$$
We have
$$|U_i| =\sum_{j \in S_i} t_j {\bf 1}_{j \notin M_i}\leq \sum_{j \in S_i} t_j\left( \sum_{j^{\prime},j^{\prime\prime}| [j,j^{\prime},j^{\prime\prime}]}t_{j^{\prime}}t_{j^{\prime\prime}} + \sum_{j^{\prime}\in A_i(j)}t_{j^{\prime}}\right).$$
Consider
$$Y= \sum_{j \in S_i} t_j  \sum_{j^{\prime}\in A_i(j)}t_{j^{\prime}}=\sum_{j \in S_i , j^{\prime}\in A_i(j)} t_{j}t_{j^{\prime}}.$$
To use Lemma 4.2 of \cite{KV2003}, we have to compute $\mathbb{E}_i(Y)$, for $i=0,1,2$. By the induction hypothesis $|S_i|\leq \sqrt{2} b_iq^{N}$ (Property \eqref{13}) and $|A_i(j)|\leq \sqrt{2}a_ib_i q^{\frac{N+1}{2}}$ (Property \eqref{12}). Then
$$\mathbb{E}_0(Y)=\mathbb{E}(Y)\leq  \sqrt{2}a_ib_i q^{\frac{N+1}{2}}\sqrt{2} b_iq^{N}p_i^2=2a_i b_i^2 \frac{\theta^2}{b_i^2q^{N+1}}q^Nq^{\frac{N+1}{2}}=2a_i\theta^2q^{\frac{N-1}{2}}.$$
To bound $\mathbb{E}_1(Y)$, notice that for a fixed $j$ there are at most $\max_{j}|A_i(j)|$ points $j^{\prime}$ such that $t_jt_{j^{\prime}}$ appears in $Y$, so
$$\mathbb{E}_1(Y)\leq  \sqrt{2}a_ib_i q^{\frac{N+1}{2}}p_i=\sqrt{2}a_i b_i q^{\frac{N+1}{2}}\frac{\theta}{b_iq^{\frac{N+1}{2}}}=\sqrt{2}a_i \theta\leq 2 i \theta^2=o(1),$$
since $a_i\simeq i\theta$.
To bound $\mathbb{E}_2(Y)$, notice that each product $t_jt_{j^{\prime}}$ can appear only at most twice in $Y$, so
$$\mathbb{E}_2(Y)\leq 2.$$
We can now apply Lemma 4.2 of \cite{KV2003} and we get that $Y$ is strongly concentrate and in particular
$$Y\leq 3a_i\theta^2q^{\frac{N-1}{2}}=o\left(\theta q^{\frac{N-1}{2}}\right)$$
with very high probability.
Now consider
$$Y=\sum_{j \in S_i} t_j \sum_{j^{\prime},j^{\prime\prime}| [j,j^{\prime},j^{\prime\prime}]}t_{j^{\prime}}t_{j^{\prime\prime}}=\sum_{j , j^{\prime},j^{\prime\prime}| [j,j^{\prime},j^{\prime\prime}]}t_jt_{j^{\prime}}t_{j^{\prime\prime}}.$$
Arguing as before, we have
\begin{eqnarray*}
\mathbb{E}_0(Y)=\mathbb{E}(Y)\leq  |S_i|^2|S_i((jj^{\prime}))|p_i^3\leq (\sqrt{2} b_iq^{N})^2 \sqrt{2}b_iq p_i^3=2\sqrt{2} b_i^3 q^{2n+1} \frac{\theta^3}{b_i^3q^{\frac{3n+3}{2}}}=2\sqrt{2}\theta^3q^{\frac{N-1}{2}}\\
\mathbb{E}_1(Y)\leq  |S_i||S_i((jj^{\prime}))|p_i^2\leq \sqrt{2} b_iq^{N} \sqrt{2}b_iq p_i^2=2 b_i^2 q^{N+1} \frac{\theta^2}{b_i^2q^{N+1}}=2\theta^2\\
\mathbb{E}_2(Y)\leq  |S_i((jj^{\prime}))|p_i\leq b_iq p_i= b_i q \frac{\theta}{b_iq^{\frac{N+1}{2}}}=o(1)\\
\mathbb{E}_3(Y)\leq  6.
\end{eqnarray*}
Then, with very high probability
$$Y\leq 4\theta^3q^{\frac{N-1}{2}}=o\left(\theta q^{\frac{N-1}{2}}\right),$$
which implies $|U_i|=o\left(\theta q^{\frac{N-1}{2}}\right)$ and therefore we have proven the bound on $|M_i|$.
\end{itemize}
\endproof

\begin{proposition}\label{Property2}
Property \eqref{2} holds.
\end{proposition}
\proof
We have to prove that
$$b_{i+1}q(1-(i+1)\log^{-13}q)\leq |S_{i+1}(\ell)|\leq b_{i+1}q(1+(i+1)\log^{-13}q).$$
Consider $L=S_i(\ell)$ and, in order to apply Corollary \ref{Corollary 4.8}, $K=14$. Let  $x \in A(S_{i}(\ell),j)$ then: $x \in \ell$,  $x \in S_{i}$, and there exists $a_1 \in A_i$ such that $[a_1xj]$ (then $x \in S_i(\ell,j)$); so, by Property \eqref{4},
$$a(L)\leq \max_{\ell^{\prime},v} |S_i(\ell^{\prime},v)|=O(ia_ib_i\sqrt{q}+i\log^{40} q).$$
By the induction hypothesis, we have that $|L|=|S_{i}(\ell)|\geq \frac{1}{2}b_i q\geq \frac{1}{2}\log^{100} q$ (recall that we are in the first phase) then, since $a_i=o(log^2 q)$,

$$a(L)\leq O(ia_ib_i\sqrt{q}+i\log^{40} q)\leq \frac{\frac{1}{2}b_i q}{\log^{38}q}\leq  \frac{|L|}{log^{2( K+5)} q}.$$

Therefore, by Corollary \ref{Corollary 4.8}, with very high probability
$$||L^{\prime}|- \mathbb{E}(|L^{\prime}|)|=||S_{i+1}(\ell)|- \mathbb{E}(|L^{\prime}|)|\leq |L|\log^{-14} q\leq 2b_iq \log^{-14} q.$$
Note that (see also Remark \ref{RemarkComp}) $\mathbb{E}(|L^{\prime}|)=\frac{b_{i+1}}{b_i}|L|$ and then, by hypothesis,
$$|\mathbb{E}(|L^{\prime}|)- b_{i+1}q|\leq ib_{i+1}q \log^{-13} q.$$
Finally we get from the two inequalities above
$$||L^{\prime}|- b_{i+1}q|=||S_{i+1}(\ell)|- b_{i+1}q|\leq ib_{i+1}q \log^{-13} q+2b_iq \log^{-14} q\leq (i+1)b_{i+1}q \log^{-13} q.$$
\endproof

\begin{proposition}\label{Property3}
Property \eqref{3} holds.
\end{proposition}
\proof
The proof is similar to the proof of Proposition \ref{Property2}.

We have to prove that
$$ |\Omega_{i+1}(\ell)|\leq b_{i+1}^{\prime}q(1+(i+1)\log^{-13}q).$$
Consider $L=\Omega_i(\ell)$ and, in order to apply Corollary \ref{Corollary 4.8}, $K=14$. Let $x \in A(\Omega_{i}(\ell),j)$ then: $x \in \ell$, $x \in \Omega_{i}$, and there exists $a_1 \in A_i$ such that $[a_1xj]$ (then $x \in \Omega_i(\ell,j)$); so, by Property \eqref{4},
$$a(L)\leq \max_{\ell^{\prime},v} |\Omega_i(\ell^{\prime},v)|=O(ia_ib_i^{\prime}\sqrt{q}+i\log^{40} q).$$
By the induction hypothesis, we have that $|L|=|\Omega_{i}(\ell)|\geq |S_{i}(\ell)|\geq \frac{1}{2}b_i q\geq \frac{1}{2}\log^{100} q$ (recall that we are in the first phase) then, since $a_i=o(log^2 q)$,
$$a(L)\leq O(ia_ib_i^{\prime}\sqrt{q}+i\log^{40} q) \leq \frac{\frac{1}{2}b_iq}{\log^{38}q}\leq \frac{|L|}{\log^{2(K+5)} q}.$$
Therefore, by Corollary \ref{Corollary 4.8}, with very high probability
$$|L^{\prime}|\leq  \mathbb{E}(|L^{\prime}|)+ |L|\log^{-14} q\leq \mathbb{E}(|L^{\prime}|)+2b_iq \log^{-14} q.$$
Note that $\mathbb{E}(|L^{\prime}|)\leq \frac{b_{i+1}^{\prime}}{b_i^{\prime}}|L|$ and then, by hypotesis,
$$|\mathbb{E}(|L^{\prime}|)- b_{i+1}q|\leq ib_{i+1}^{\prime}q \log^{-13} q.$$
Finally we get from the two inequalities above
$$|L^{\prime}|=|\Omega_{i+1}(\ell)|\leq b_{i+1}^{\prime}q+ ib_{i+1}^{\prime}q \log^{-13} q+2b_i^{\prime}q \log^{-14} q\leq b_{i+1}^{\prime}q +(i+1)b_{i+1}^{\prime}q \log^{-13} q.$$
\endproof

\begin{proposition}\label{Property4}
Property \eqref{4} holds.
\end{proposition}
\proof
We have to show that for every $j$
$$|S_{j}(\ell,v)|\leq 8ja_jb_j\sqrt{q}+j\log^{40}q.$$
Consider $L=S_{i}(\ell,v)$. By hypotesis
$$|S_{i}(\ell,v)|\leq 8ia_ib_i\sqrt{q}+i\log^{40}q.$$
Notice that if $x\in S_{i+1}(\ell,v)$ but not in $L^{\prime}=L\cap S_{i+1}$, it means that in the line $(xv)$ some point is chosen and belongs to $A_{i+1}$; also,
$$|S_{i+1}(\ell,v)|= |L^{\prime}|+|S_{i+1}(\ell,v)\setminus L^{\prime}|.$$
Recall that $$A(L,j)=\{t \in L | (tj) \cap A_i\neq \emptyset\}$$ and
$$a(L)=\max \left\{ \max_{j \in \Omega_i} |A(L,j)|,|L|\log^{-100}q\right\}.$$
Let  $x \in A(S_{i}(\ell,v),j)$, then: $x \in \ell$, $x \in S_{i}$, there exists $a_1 \in A_i$ such that $[a_1xv]$, and there exists $a_2 \in A_i$ such that $[a_2xj]$. So, $x \in S_i(\ell,v,j)$. Therefore
$$a(L)\leq \max_{\ell,u,v} |S_i(\ell,v,u)|\leq i\log^{4} q\leq \log^{7}q.$$
Observe we used  $i\leq \log^3 q$, which is true by hypothesis. Moreover (see Remark \ref{RemarkComp})
$$\mathbb{E}(|L^{\prime}|)=|L|(1-P_i^u)\leq 8ia_ib_{i+1}\sqrt{q}+i\log^{40}q,$$
since $b_{i+1}=b_i(1-P_i^u)$. By Lemma \ref{Lemma 4.7}, with very high probability
$$|L^{\prime}|\leq \mathbb{E}(L^{\prime})+a(L)^{1/2}|L|^{1/2}\log^5 q\leq 8ia_ib_{i+1}\sqrt{q}+i\log^{40}q+|L|^{1/2}\log^{9}q.$$
For the other set, since in the line $(xv)$ some point $j$ has $t_j=1$, we have
$$|S_{i+1}(\ell,v)\setminus L^{\prime}|\leq \sum _{x \in S_i(\ell)} \sum_{j \in (xv)} t_j=Y.$$
As before,
\begin{eqnarray*}
\mathbb{E}_0(Y)=\mathbb{E}(Y)\leq  p_i|S_i(\ell^{\prime})|^2\leq (\sqrt{2} b_iq)^2 p_i=2 b_i^2 q^{2} \frac{\theta}{b_iq^{\frac{N+1}{2}}}=2b_i\theta q^{-\frac{N-3}{2}}=o(1)\\
\mathbb{E}_1(Y)\leq  1.
\end{eqnarray*}
(recall that $N\geq 3$). So, with very high probability
$$|S_{i+1}(\ell,v)\setminus L^{\prime}|=O(\log^3 q).$$
Therefore
$$(8ia_ib_{i+1}\sqrt{q}+i\log^{40}q)^{1/2}\log^{9}q\leq 8a_ib_{i+1}\sqrt{q},
$$
so
$$(8ia_ib_{i+1}\sqrt{q}+i\log^{40}q)^{1/2}\log^{9}q+\log^4 q\leq 8a_ib_{i+1}\sqrt{q}+\log^{40} q;$$
then
$$|S_{i+1}(\ell,v)|\leq 8ia_ib_{i+1}\sqrt{q}+i\log^{40}q+ 8a_ib_{i+1}\sqrt{q}+\log^{40} \leq$$
$$\leq 8(i+1)a_{i}b_{i+1}\sqrt{q}+(i+1)\log^{40}q\leq 8(i+1)a_{i+1}b_{i+1}\sqrt{q}+(i+1)\log^{40}q,$$
since $a_{i+1}\simeq a_{i} +\theta$.
\endproof

\begin{proposition}\label{Property5}
Property \eqref{5} holds.
\end{proposition}
\proof
We have to prove that $|S_{i+1}(\ell,u,v)|\leq (i+1)\log^4 q$. Recall that a point $x$ is in $S_{i+1}(\ell,u,v)$ if it lies on the line $\ell$ and there exist two points $a_1,a_2\in A_{i+1}$ such that $[x,a_1,u]$ and $[x,a_2,v]$. Let $B_i(\ell,u)$ be the set of all the $x \in \ell$ such that there exists $b\in B_i$ with $[x,b,u]$ and $B_i(\ell,u,v)=B_i(\ell,u)\cap B_i(\ell,v)$. If $x \in S_{i+1}(\ell,u,v) \setminus S_{i}(\ell,u,v)$ it means that
$$x\in B_i(\ell,u,v) \cup (S_i(\ell,u)\cap B_i(\ell,v)) \cup (S_i(\ell,v)\cap B_i(\ell,u)),$$
so
$$|S_{i+1}(\ell,u,v)|\leq |S_{i}(\ell,u,v)|+|B_i(\ell,u,v)|+ |S_i(\ell,u)\cap B_i(\ell,v)| + |S_i(\ell,v)\cap B_i(\ell,u)|.$$
Let $x \in B_i(\ell,u,v)$, then $x \in S_i(\ell)$ and there exist $j,j^{\prime}\in B_i$ with $[x,u,j]$ and $[x,v,j^{\prime}]$ which are ``chosen'', that is $t_j=t_{j^{\prime}}=1$. Then
$$|B_i(\ell,u,v)|\leq \sum_{x \in S_i(\ell)} \left( \sum_{j | [xuj]}t_j\right) \left( \sum_{j^{\prime} | [xvj^{\prime}]}t_j^{\prime}\right)=Y$$
and 
\begin{eqnarray*}
\mathbb{E}_0(Y)=\mathbb{E}(Y)\leq  |S_i(\ell^{\prime})|^3 p_i^2\leq (\sqrt{2} b_iq)^3  p_i^2\leq 4 b_i^3 q^{3} \frac{\theta^2}{b_i^2q^{N+1}}=4b_i \theta^2q^{-N+2}=o(1)\\
\mathbb{E}_1(Y)\leq  |S_i(\ell^{\prime})|p_i\leq \sqrt{2} b_iq p_i=\sqrt{2} b_i q \frac{\theta}{b_i q^{\frac{N+1}{2}}}=\sqrt{2}\theta q^{-\frac{N-1}{2}}=o(1)\\
\mathbb{E}_2(Y)\leq  2.
\end{eqnarray*}
Thus, with very high probability (as noted after   \cite[Corollary 4.3]{KV2003}) $|B_i(\ell,u,v)|=O(\log^3 q)$.

For the second term,
$$|S_i(\ell,u)\cap B_i(\ell,v)|\leq \sum_{x \in S_i(\ell)} \left( \sum_{j\in A_i(u) | [xuj]}t_j\right) \left( \sum_{j^{\prime} | [xvj^{\prime}]}t_j^{\prime}\right)=Y$$
and 
\begin{eqnarray*}
\mathbb{E}_0(Y)=\mathbb{E}(Y)\leq  |S_i(\ell^{\prime})|^3p_i^2\leq 4 b_i^3q^3 p_i^2=4 b_i^3 q^{3} \frac{\theta^2}{b_i^2q^{N+1}}=4b_i \theta^2q^{-N+2}=o(1)\\
\mathbb{E}_1(Y)\leq  |S_i(\ell)|p_i=o(1)\\
\mathbb{E}_2(Y)\leq  2.
\end{eqnarray*}
As above, $|S_i(\ell,u)\cap B_i(\ell,v)|=O(\log^3 q)$. The same for $|S_i(\ell,v)\cap B_i(\ell,u)|$. Since by induction $|S_{i}(\ell,u,v)|\leq i\log^4 q$, then
$$|S_{i+1}(\ell,u,v)|\leq (i+1)\log^4q.$$
\endproof

\begin{proposition}\label{Property6}
Property \eqref{6} holds.
\end{proposition}
\proof
We have to prove that
$$|A_{i+1}(u,v)|\leq (i+1)b_{i+1}q+(i+1)\log^{40}q.$$
A point in $A_{i+1}(u,v)$ can be contained in $A_i(u,v)$, $B_i(u,v)$, $A_{i}(u) \cap B_i (v)$ or $A_{i}(v) \cap B_i (u)$. Then
$$|A_{i+1}(u,v)|\leq |A_i(u,v)|+ |B_i(u,v)| + |A_{i}(u) \cap B_i (v)|+|A_{i}(v) \cap B_i (u)|.$$

\begin{itemize}
\item By the induction hypotesis $|A_{i}(u,v)|\leq ib_iq+i\log^{40} q$.
\item $B_i(u,v)$: for any $x \in B_i(u,v)$ there exist $j,j^{\prime} \in B_i$ such that $[xuj]$, $[xvj^{\prime}]$. So,
$$|B_i(u,v)|\leq \sum_{x \in S_i}\left(\sum_{j |[xuj]} t_j\right)\left(\sum_{j^{\prime}|[xvj^{\prime}]} t_{j^{\prime}}\right)=Y.$$
In order to evaluate $\mathbb{E}_i(Y)$, notice that if $j$ is fixed, then $x\in S_i((ju))$, while $j^{\prime}$ belongs to  $S_i((xv)).$
Then
\begin{eqnarray*}
\mathbb{E}_0(Y)=\mathbb{E}(Y)\leq  p_i^2|S_i| (\max_{\ell} |S_i(\ell)|)^2=\frac{\theta^2}{b_i^2 q^{N+1}} 2 b_i q^N(2b_i q)^2 =8 b_i \theta^2 q\\
\mathbb{E}_1(Y)\leq   p_i (\max_{\ell} |S_i(\ell)|)^2=\frac{\theta}{b_i q^{\frac{N+1}{2}}}(2b_i q)^2 =4 b_i \theta q^{-\frac{N-3}{2}}=o(1)\\
\mathbb{E}_2(Y)\leq 2.
\end{eqnarray*}
By Lemma \ref{Lemma 4.2} and using $\lambda=\log^{3/2}q $ we get $|B_i(X)|\leq 8 b_i \theta^2 q + 4 \theta \sqrt{q}\log^3 q=o(b_i q\log^{-3}q)$.
\item $A_i(u)\cap B_i(v)$: for any $x \in A_i(u)\cap B_i(v)$ there exist $j_1\in A_i$ and $j_{2} \in B_i$ such that $[xuj_1]$ and $[xvj_2]$.
$$|A_i(u)\cap B_i(v)|\leq \sum_{x \in (uj_1)}\left(\sum_{j_{2}|[xvj_{2}]} t_{j_{2}}\right)=Y.$$
Then
\begin{eqnarray*}
\mathbb{E}_0(Y)\leq  p_i|A_i(u)| (\max_{\ell} |S_i(\ell)|)=\frac{\theta}{b_i q^{\frac{N+1}{2}}} 2 i a_ib_iq^{\frac{N+1}{2}}(2b_i q) =4 a_i b_i \theta q\\
\mathbb{E}_1(Y)=1.
\end{eqnarray*}
By Lemma \ref{Lemma 4.2} and using $\lambda=\log^{1/2}q $ we get $|A_i(u)\cap B_i(v)|\leq 4 b_i \theta^2 q + 2 i\theta \sqrt{q}=o(b_i q\log^{-3}q)$.
Therefore
$$|A_{i+1}(u,v)|\leq ib_iq+i\log^{40} q+o(b_iq\log^{-3}q)\leq  (i+1)b_{i+1}q+(i+1)\log^{40} q.$$
\end{itemize}
\endproof

\begin{proposition}\label{Property9}
Property \eqref{9} holds.
\end{proposition}
\proof
The proof is very similar to the proof of Proposition \ref{Property4}.

We have to show that for every $j$
$$|\Omega_{j}(\ell,v)|\leq 8ja_jb_j^{\prime}\sqrt{q}+j\log^{40}q.$$
Consider $L=\Omega_{i}(\ell,v)$. By hypotesis
$|\Omega_{i}(\ell,v)|\leq 8ia_ib_i^{\prime}\sqrt{q}+i\log^{40}q.$
Notice that if $x\in \Omega_{i+1}(\ell,v)$ but not in $L^{\prime}=L\cap \Omega_{i+1}$, it means that in the line $(xv)$ some point is chosen and belongs to $A_{i+1}$; also,
$$|\Omega_{i+1}(\ell,v)|= |L^{\prime}|+|\Omega_{i+1}(\ell,v)\setminus L^{\prime}|.$$
Recall that $$A(L,j)=\{t \in L | (tj) \cap A_i\neq \emptyset\}$$ and
$$a(L)=\max \left\{ \max_{j \in \Omega_i} |A(L,j)|,|L|\log^{-100}q\right\}.$$
Let $x \in A(\Omega_{i}(\ell,v),j)$, then: $x \in \ell$, $x \in \Omega_{i}$, there exists $a_1 \in A_i$ such that $[a_1xv]$, and there exists $a_2 \in A_i$ such that $[a_2xj]$. So, $x \in \Omega_i(\ell,v,j)$. Therefore
$$a(L)\leq \max_{\ell,u,v} |\Omega_i(\ell,v,u)|\leq i\log^{4} q\leq \log^{7}q.$$
Notice that we used that $i\leq \log^3 q$ which is true by hypothesis. Moreover
$$\mathbb{E}(|L^{\prime}|)\leq|L|(1-P_i^{\ell})\leq 8ia_ib_{i+1}^{\prime}\sqrt{q}+i\log^{40}q,$$
since $b_{i+1}^{\prime}=b_i^{\prime}(1-P_i^{\ell})$. By Lemma \ref{Lemma 4.7}, with very high probability
$$|L^{\prime}|\leq \mathbb{E}(L^{\prime})+a(L)^{1/2}|L|^{1/2}\log^5 q\leq 8ia_ib_{i+1}^{\prime}\sqrt{q}+i\log^{40}q+|L|^{1/2}\log^{9}q.$$
For the other set, since in the line $(xv)$ some point $j$ has $t_j=1$, we have
$$|\Omega_{i+1}(\ell,v)\setminus L^{\prime}|\leq \sum _{x \in \ell} \sum_{j \in (xv)} t_j=Y.$$
As before
\begin{eqnarray*}
\mathbb{E}_0(Y)=\mathbb{E}(Y)\leq  p_i|\Omega_i(\ell^{\prime})|^2\leq (\sqrt{2} b_i^{\prime}q)^2 p_i=2 (b_i^{\prime})^2 q^{2} \frac{\theta}{b_iq^{\frac{N+1}{2}}}\leq 4b_i\theta q^{-\frac{N-3}{2}}=o(1)\\
\mathbb{E}_1(Y)\leq  1.
\end{eqnarray*}
(recall that $N\geq 3$). So, with very high probability
$$|\Omega_{i+1}(\ell,v)\setminus L^{\prime}|=O(\log^3 q).$$
Therefore
$$(8ia_ib_{i}^{\prime}\sqrt{q}+i\log^{40}q)^{1/2}\log^{9}q\leq 8a_ib_{i+1}^{\prime}\sqrt{q},
$$
so
$$(8ia_ib_{i}^{\prime}\sqrt{q}+i\log^{40}q)^{1/2}\log^{9}q+\log^{40} q\leq 8a_ib_{i+1}^{\prime}\sqrt{q}+\log^{40} q;$$
then
$$|\Omega_{i+1}(\ell,v)|\leq 8ia_ib_{i+1}^{\prime}\sqrt{q}+i\log^{40}q+ 8a_ib_{i+1}^{\prime}\sqrt{q}+\log^{40}q =$$
$$= 8(i+1)a_{i}b_{i+1}^{\prime}\sqrt{q}+(i+1)\log^{40}q\leq 8(i+1)a_{i+1}b_{i+1}^{\prime}\sqrt{q}+(i+1)\log^{40}q,$$
since $a_{i+1}\simeq a_{i} +\theta$.
\endproof

\begin{proposition}\label{Property10}
Property \eqref{10} holds.
\end{proposition}

\proof
The proof is similar to the proof of Proposition \ref{Property5}.

We have to prove that $|\Omega_{i+1}(\ell,u,v)|\leq (i+1)\log^4 q$. Recall that a point $x\in \Omega_i$ is in $\Omega_{i+1}(\ell,u,v)$ if it lies on the line $\ell$ and there exist two points $a_1,a_2\in A_{i+1}$ such that $[x,a_1,u]$ and $[x,a_2,v]$. Let $B_i(\ell,u)$ be the set of all the $x \in \ell$ such that there exists $b\in B_i$ with $[x,b,u]$ and $B_i(\ell,u,v)=B_i(\ell,u)\cap B_i(\ell,v)$. If $x \in \Omega_{i+1}(\ell,u,v) \setminus \Omega_{i}(\ell,u,v)$ it means that
$$x\in B_i(\ell,u,v) \cup (\Omega_i(\ell,u)\cap B_i(\ell,v)) \cup (\Omega_i(\ell,v)\cap B_i(\ell,u)),$$
so
$$|\Omega_{i+1}(\ell,u,v)|\leq |\Omega_{i}(\ell,u,v)|+|B_i(\ell,u,v)|+ |\Omega_i(\ell,u)\cap B_i(\ell,v)| + |\Omega_i(\ell,v)\cap B_i(\ell,u)|.$$
Let  $x \in B_i(\ell,u,v)$, then there exist $j,j^{\prime}\in B_i$ with $[x,u,j]$ and $[x,v,j^{\prime}]$ which are ``chosen'', that is $t_j=t_{j^{\prime}}=1$. Then
$$|B_i(\ell,u,v)|\leq \sum_{x \in \ell} \left( \sum_{j | [xuj]}t_j\right) \left( \sum_{j^{\prime} | [xvj^{\prime}]}t_j^{\prime}\right)=Y.$$
We already know from the proof of the Main Theorem, that
$$\log \frac{b_i^{\prime}}{b_i}=\log\prod_{j=0}^{i}\frac{1-P_j^{\ell}}{1-P_j^{u}}=O(i\log^{-3}q)=o(1) \Longrightarrow b_i^{\prime}\leq 2 b_i.$$
Then
\begin{eqnarray*}
\mathbb{E}_0(Y)=\mathbb{E}(Y)\leq  |\Omega_i(\ell^{\prime})|^3 p_i^2\leq (\sqrt{2} b_i^{\prime}q)^3  p_i^2\leq 4 (b_i^{\prime})^3 q^{3} \frac{\theta^2}{b_i^2q^{N+1}}\leq 32 b_i \theta^2q^{-N+2}=o(1)\\
\mathbb{E}_1(Y)\leq  |\Omega_i(\ell^{\prime})|p_i\leq \sqrt{2} b_i^{\prime}q p_i=\sqrt{2} b_i^{\prime} q \frac{\theta}{b_i q^{\frac{N+1}{2}}}=2\sqrt{2}\theta q^{-\frac{N-1}{2}}=o(1)\\
\mathbb{E}_2(Y)\leq  2.
\end{eqnarray*}
With very high probability (as noted after \cite[Corollary 4.3]{KV2003}) $|B_i(\ell,u,v)|=O(\log^3 q)$. Also, consider 

$$|\Omega_i(\ell,u)\cap B_i(\ell,v)|\leq \sum_{x \in \ell} \left( \sum_{j\in A_i(u) | [xuj]}t_j\right) \left( \sum_{j^{\prime} | [xvj^{\prime}]}t_j^{\prime}\right)=Y.$$
Then
\begin{eqnarray*}
\mathbb{E}_0(Y)=\mathbb{E}(Y)\leq  |\Omega_i(\ell^{\prime})|^3 p_i^2\leq (\sqrt{2} b_i^{\prime}q)^3 p_i^2\leq 4 (b_i^{\prime})^3 q^3 \frac{\theta^2}{b_i^2q^{N+1}}\leq 32 b_i \theta^2q^{-N+2}=o(1)\\
\mathbb{E}_1(Y)\leq  |\Omega_i(\ell)|p_i=o(1)\\
\mathbb{E}_2(Y)\leq  2.
\end{eqnarray*}
As above, $|\Omega_i(\ell,u)\cap B_i(\ell,v)|=O(\log^3 q)$. The same for $|\Omega_i(\ell,v)\cap B_i(\ell,u)|$. Since by induction $|\Omega_{i}(\ell,u,v)|\leq i\log^4 q$, thus
$$|\Omega_{i+1}(\ell,u,v)|\leq (i+1)\log^4q.$$
\endproof

\section{Second phase}\label{Phase2}

\begin{remark}\label{ProprValide}
Property \ref{Property6} is still valid in Phase two, since in the proof we used properties which are still valid in this phase. Also, let $\bar i$ such that
$b_{\bar i+1} <\frac{\log ^{c_1} q}{q}\leq b_{\bar i}$. Then for $i> \bar i+1$
$$|S_{i}(\ell)| \leq |S_{\bar i+1}(\ell)| \leq  2b_{\bar i +1}q\leq 2\log^{c_1} q$$
and
$$|\Omega_{i}(\ell)| \leq |\Omega_{\bar i+1}(\ell)| \leq  2b^{\prime}_{\bar i +1}q\leq  4b_{\bar i +1}q\leq \log^{c_1} q.$$
\end{remark}

\begin{proposition}\label{Property1bis}
Property \eqref{1bis} holds.
\end{proposition}
\proof
The proof is the same of Proposition \ref{Property1}.
\endproof

\begin{proposition}\label{Property2bis}
Property \eqref{2bis} holds.
\end{proposition}
\proof
We have to prove that with very high probability
$$\frac{1}{2}b_{i+1}^2q^{N+1}(1-3(i+1)\log^{-13}q)\leq |T_{i+1}(v)|\leq \frac{1}{2} b_{i+1}^2q^{N+1}(1+3(i+1)\log^{-13}q).$$
Recall that $T_i(v)=\{\{x,y\}| x,y \in S_i \setminus\{v\}, [xyv]\},$
then
$$\mathbb{E}(T_{i+1}(v))=\sum_{x,y \in S_i, [xyv]} Pr(x,y \in S_{i+1}).$$
Therefore, by Lemma \ref{Lemma 5.3},
$$|T_i(v)|Pr(x\in S_{i+1}|x \in S_i)^2(1-o(\log^{-13}q))\leq \mathbb{E}(T_{i+1}(v)) \leq$$
$$\leq |T_i(v)|Pr(x\in S_{i+1}|x \in S_i)^2(1+o(\log^{-13}q)).$$
By Remark \ref{RemarkComp} $Pr(x \in S_{i+1}|x \in S_i)=1-P_i^u$, so $b_iPr(x \in S_{i+1}|x \in S_i)=b_{i+1}$, and by the induction hypotesis
$$\frac{1}{2} q^{N+1}b_{i+1}^2(1-3i\log^{-13}q)(1-o(\log^{-13}q))\leq \mathbb{E}(T_{i+1}(v)) \leq$$
$$\leq\frac{1}{2} q^{N+1}b_{i+1}^2(1+3i\log^{-13}q)(1+o(\log^{-13}q)).$$
Consider
$$T_i^{\prime}(v) =\{ x^{m(x)}| x \in S_i \setminus \{v\}\},$$
where $m(x)=|S_i((x,v))|-2$. It is easy to see that $|T_i^{\prime}(v)|=2|T_i(v)|$. Let $L=T_i^{\prime}(v)$ and $L^{\prime}=T_{i+1}^{\prime}(v)$, it holds that 

$$|A(L,u)|=|\{t \in L | (tu)\cap A \neq \emptyset\}|\leq |A_i(u)|\max_x m(x),$$
since every element $x\in A_i(u)$ is repeated in $A(L,u)$ exactly $m(x)$ times. Moreover, since we are in the second phase and $b_iq\leq\log^{c_1}q$,
$$|A_i(u)|\max_x m(x)\leq|A_i(u)|\max_{\ell} |S_i(\ell)|\leq \sqrt{2}a_ib_iq^{\frac{N+1}{2}}\sqrt{2}\log^{c_1}q=2a_ib_iq^{\frac{N+1}{2}}\log^{c_1}q.$$
Since $|L|=|T_i^{\prime}(v)|\geq \frac{1}{3}b_i^2q^{N+1}\geq \frac{a(L)}{\log^{100}q}$, by Lemma \ref{Lemma 4.7} we have that with very high probability
$$\left||L^{\prime}|-\mathbb{E}(|L^{\prime}|)\right|\leq (|L|a(L))^{1/2} \log^5 q\leq \left(2b_i^2q^{N+1}2a_ib_iq^{\frac{N+1}{2}}\log^{c_1}q\right)^{1/2}\log^5 q=$$
$$=2a_i^{1/2}b_i^{3/2}q^{\frac{3(n+1)}{4}}\log ^{55}q\leq b_i^{3/2}q^{\frac{3(n+1)}{4}}\log ^{60}q.$$
Since $b_iq^{\frac{N+1}{2}}\geq \log^{c} q$ and $2b_{i+1}\geq b_i$, we have that
$$\frac{b_i^{3/2}q^{\frac{3(n+1)}{4}}\log ^{60}q}{ b_{i+1}^2 q^{N+1}\log^{-13}q}\leq 4\frac{b_i^{3/2}q^{\frac{3(n+1)}{4}}\log ^{60}q}{ b_{i}^2 q^{N+1}\log^{-13}q}\leq 4\left(b_i q^{\frac{N+1}{2}}\right)^{-1/2}\log ^{73}q\leq4 \log^{-77} q.$$
So $b_i^{3/2}q^{\frac{3(n+1)}{4}}\log ^{60}q=o( b_{i+1}^2 q^{N+1}\log^{-13}q)$.
Therefore with very high probability
$$b_{i+1}^2 q^{N+1}(1-3i\log^{-13} q)(1-o(\log^{-13}q))-o( b_{i+1}^2 q^{N+1}\log^{-13}q) \leq |L^{\prime}|\leq$$
$$\leq b_{i+1}^2 q^{N+1}(1+3i\log^{-13} q)(1+o(\log^{-13}q))+o( b_{i+1}^2 q^{N+1}\log^{-13}q).$$
We have that
$$b_{i+1}^2 q^{N+1}(1-3i\log^{-13} q)(1-o(\log^{-13}q))-o( b_{i+1}^2 q^{N+1}\log^{-13}q)=$$
$$=b_{i+1}^2 q^{N+1}(1-3i\log^{-13} q-3o(\log^{-13}q))\geq b_{i+1}^2 q^{N+1}(1-3(i+1)\log^{-13} q)$$
and
$$b_{i+1}^2 q^{N+1}(1+3i\log^{-13} q)(1+o(\log^{-13}q))+o( b_{i+1}^2 q^{N+1}\log^{-13}q)=$$
$$=b_{i+1}^2 q^{N+1}(1+3i\log^{-13} q+3(\log^{-13}q))\leq b_{i+1}^2 q^{N+1}(1+3(i+1)\log^{-13} q).$$
Recalling that $|L^{\prime}|=2|T_i(v)|$ we have the result.
\endproof

\begin{proposition}\label{Property3bis}
Property \eqref{3bis} holds.
\end{proposition}
\proof
We have to prove that
$$a_{i+1}b_{i+1}q^{\frac{N+1}{2}}(1-O(i\theta^2) )\leq |A_{i+1}(v)|\leq a_{i+1}b_{i+1}q^{\frac{N+1}{2}}(1+O(i\theta^2) ).$$
Let $U_i=B_i \setminus M_i$ and $U_i(v),B_i(v),M_i(v)$ the set of points $x \in S_i$ such that there exists $u\in U_i,B_i,M_i$ respectively with $[xuv]$, and $B_i^{\prime}(v)=B_i(v)\cap S_{i+1}$, $M_i^{\prime}(v)=M_i(v)\cap S_{i+1}$. Let $L=A_i(v)$, $L^{\prime}=L\cap S_{i+1}$, $v \in \Omega_{i+1}$, then
$$|A_{i+1}(v)|=|L^{\prime}|+|M_i^{\prime}(v)|.$$
In the following we will give some estimations for the sizes of $B_{i}(v)$, $U_{i}(v)$, and $B_{i}^{\prime}(v)$.

\begin{itemize}
\item {\bf $|B_{i}(v)|\leq \theta b_iq^{\frac{N+1}{2}}(1 + o(\log^{-10}q)).$}\\
If $x \in B_i(v)$, then there exists $j\in B_i$ such that $[jxv]$ and $t_j=1$. Then
$$|B_i(v)|\leq \sum_{x \in S_i} \left( \sum_{j\in S_i | [jxv]} t_j\right)=\sum_{j \in S_i}t_j m(j)=Y,$$
where $m(j)=|S_i((vj))\setminus\{v,j\}|$. Notice that $\sum_{j \in S_i} m(j)=|T_{i}^{\prime}(v)|=2|T_i(v)|$. By Property \eqref{2bis},
\begin{eqnarray*}
\mathbb{E}_0(Y)=\mathbb{E}(Y)=  2p_i|T_i(v)|\geq p_ib_i^2q^{N+1}(i-o(\log^{-10}q))\simeq \theta b_iq^{\frac{N+1}{2}}\geq \log^{150}q\\
\mathbb{E}_1(Y)\leq   \max_x m(x)\leq \max_{\ell}|S_i(\ell)|\leq 2\log^{c_1}q
\end{eqnarray*}
By Lemma \ref{Lemma 4.2}, with very high probability,
$$Y=|B_i(v)|\leq \mathbb{E}(Y)+ (2\mathbb{E}(Y)\log^{c_1}q )^{1/2}\log^2 q\leq\theta b_i q^{\frac{N+1}{2}}(1+o(\log^{-10}q)).$$

\item  {\bf $|U_{i}(v)|\leq 6\theta^2 a_ib_iq^{\frac{N+1}{2}}.$}\\
Note that for any $x\in U_i(v)$ there exists $j\in (xv)$ different from $x,v$ such that $j \in B_i\setminus M_i$; so
$$|U_i(v)| \leq \sum_{x\in S_i | |(xv)|\geq 3} \quad \sum _{j\in S_i | [jxv]} t_j {\bf 1}_{\{j \notin M_i\}}.$$
Moreover, from the description of the algorithm $j\notin M_i$ if and only if  $[jj^{\prime},j^{\prime\prime}]$ for either some $j^{\prime},j^{\prime\prime}\in B_i$ or $j^{\prime} \in B_i$ and $j^{\prime\prime}\in A_i$. Then
$$ \sum_{x\in S_i | |(xv)|\geq 3} \quad \sum _{j\in S_i | [jxv]} t_j {\bf 1}_{\{j \notin M_i\}}\leq  \sum_{x \in S_i | |(xv)|\geq 3} \quad \sum _{j| [jxv]} t_j \left( \sum _{j^{\prime},j^{\prime\prime}| [jj^{\prime}j^{\prime\prime}]} t_{j^{\prime}}t_{j^{\prime\prime}}+\sum _{j^{\prime}\in A_i(j)}t_{j^{\prime}}\right).$$
Consider
\begin{equation}
Y_1=\sum_{x\in S_i | |(xv)|\geq 3} \; \sum _{j\in S_i | [jxv]} t_j  \; \sum _{j^{\prime},j^{\prime\prime}|[jj^{\prime}j^{\prime\prime}]} t_{j^{\prime}}t_{j^{\prime\prime}}$$
and
$$Y_2= \sum_{x\in S_i | |(xv)|\geq 3} \; \sum _{j\in S_i | [jxv]} t_j \; \sum _{j^{\prime}\in A_i(j)}t_{j^{\prime}}.
\end{equation}
We have that
\begin{eqnarray*}
\mathbb{E}_0(Y_1)\simeq p_i^3|S_i|\max_{\ell}|S_i(\ell)|\max_j |T_i(j)|\simeq \left(\frac{\theta^3}{b_i^3q^{\frac{3N+3}{2}}}\right) b_i^4q^{2N+2}\simeq \theta^3b_i q^{\frac{N+1}{2}}\\
\mathbb{E}_0(Y_2)\simeq p_i^2|S_i|\max_x |A_i(x)|\max_{\ell}|S_i(\ell)|\simeq \left(\frac{\theta^2}{b_i^2q^{N+1}}\right) a_ib_i^3q^{\frac{3N+3}{2}}\simeq a_i b_i \theta^2q^{\frac{N+1}{2}}.\\
\end{eqnarray*}
It is easy to prove that with very high probability
$$Y_1\leq 5\theta^3b_iq^{\frac{N+1}{2}}, \qquad Y_2\leq 5\theta^2a_ib_iq^{\frac{N+1}{2}};$$
therefore
$$|U_{i}(v)|\leq 6\theta^2 a_ib_iq^{\frac{N+1}{2}}.$$

\item {\bf $|B_{i}^{\prime}(v)|\geq \theta b_iq^{\frac{N+1}{2}}(1 - o(\log^{-10}q)) - 14a_i\theta^2b_i q^{\frac{N+1}{2}}.$}\\
A point $x \in B_i^{\prime}(v)$ if and only if $x \in S_{i+1}$ and there exists at least one point in the line $(xv)$ different from $x$ and $v$ belonging to $B_i$. Let $B_i^{\prime\prime}(v)=\{x \in S_{i+1}| \exists ! j \in B_{i}, j\neq x,v | [jxv]\}$. Then
$$|B_i^{\prime}(v)|\geq \sum _{x \in S_i| |(xv)|\geq 3} {\bf 1}_{\{x \in S_{i+1}\}} \sum_{j| [jxv]} t_j \left(1 -\sum_{j^{\prime}\neq j | [j^{\prime} xv]}t_{j^{\prime}}\right)$$
and 
$$|B_i^{\prime}(v)|\geq -|U_i(v)|+\sum _{x \in S_i| |(xv)|\geq 2} \left(1-\left(\sum_{g \in A_i(x)} t_g +\sum_{g,g^{\prime}|[gg^{\prime}x]}t_g t_{g^{\prime}}\right)\right) \sum_{j | [jxv]} t_j \left(1 -\sum_{j^{\prime}\neq j| [j^{\prime} xv]}t_{j^{\prime}}\right) .$$
In fact $\sum_{j | [jxv]} t_j \left(1 -\sum_{j^{\prime}\neq j| [j^{\prime} xv]}t_{j^{\prime}}\right)$ is different from zero only if $|(xv)|= 3$; in this case $$\left(1-\left(\sum_{g \in A_i(x)} t_g +\sum_{g,g^{\prime}|[gg^{\prime}x]}t_g t_{g^{\prime}}\right)\right)$$ is different from zero if and only if  $x \in M_i(v)\cap U_i(v)$.


Consider
$$Y= \sum _{x \in S_i | |(xv)|\geq 2}\quad \sum_{j | [jxv]} t_j.$$
Using the same calculations as above, we get that with very high probability
$$|Y|\geq \mathbb{E}_0(Y)(1-o(\log^{-10}q))\geq \theta b_iq^{\frac{N+1}{2}}(1-o(\log^{-10}q))(1-o(\log^{-10}q))\geq$$
$$\geq \theta b_iq^{\frac{N+1}{2}}(1-o(\log^{-10}q)).$$
Consider the following
\begin{eqnarray*}
Y_1= -\sum _{x \in S_i| |(xv)|\geq 2}\quad \sum_{j | [jxv]}\quad \sum_{j^{\prime}\neq j | [j^{\prime} xv]} t_jt_{j^{\prime}},\\
Y_2= -\sum _{x \in S_i| |(xv)|\geq 2}\quad \sum_{g \in A_i(x)}\quad\sum_{j\in S_i | [jxv]} t_gt_{j},\\
Y_3= \sum _{x \in S_i| |(xv)|\geq 2}\quad \sum_{g \in A_i(x)}\quad\sum_{j\in S_i | [jxv]} \quad\sum_{j^{\prime}\neq j | [j^{\prime} xv]} t_gt_{j}t_{j^{\prime}},\\
Y_4= -\sum _{x \in S_i| |(xv)|\geq 2}\quad \sum_{g,g^{\prime}| [g g^{\prime} x]}\quad\sum_{j\in S_i | [jxv]} t_gt_{g^{\prime}}t_{j},\\
Y_5= \sum _{x \in S_i| |(xv)|\geq 2}\quad \sum_{g,g^{\prime}| [g g^{\prime} x]}\quad\sum_{j\in S_i | [jxv]} \quad\sum_{j^{\prime}\neq j | [j^{\prime} xv]} t_gt_{g^{\prime}}t_{j}t_{j^{\prime}}.
\end{eqnarray*}
We have that
\begin{eqnarray*}
\mathbb{E}_0(Y_1)\simeq p_i^2|S_i|\max_{\ell}|S_i(\ell)|\simeq \left(\frac{\theta^2}{b_i^2q^{N+1}}\right) b_i^3q^{N+2}\simeq \theta^2b_i q,\\
\mathbb{E}_0(Y_2)\simeq p_i^2|S_i|\max_x |A_i(x)|\max_{\ell}|S_i(\ell)|\simeq \left(\frac{\theta^2}{b_i^2q^{N+1}}\right) a_ib_i^3q^{\frac{3n+3}{2}}\simeq a_i b_i \theta^2q^{\frac{N+1}{2}},\\
\mathbb{E}_0(Y_3)\simeq p_i^3|S_i|\max_x |A_i(x)|\left(\max_{\ell}|S_i(\ell)|\right)^2\simeq \left(\frac{\theta^3}{b_i^3q^{\frac{3n+3}{2}}}\right)a_ib_i^4q^{\frac{3n+5}{2}}\simeq a_ib_i\theta^3 q,\\
\mathbb{E}_0(Y_4)\simeq p_i^3|S_i|\left(\max_{\ell}|S_i(\ell)|\right)^3\simeq \left(\frac{\theta^3}{b_i^3q^{\frac{3n+3}{2}}}\right)b_i^4q^{N+3}=o(1),\\
\mathbb{E}_0(Y_5)\simeq p_i^4|S_i|\left(\max_{\ell}|S_i(\ell)|\right)^4\simeq \left(\frac{\theta^4}{b_i^4q^{2n+2}}\right)b_i^4q^{N+4}=o(1).\\
\end{eqnarray*}
Then it is easy to see that the maximum of the expectations is $\mathbb{E}_0(Y_2)$ and the maximum of the errors is $\simeq a_i b_i \theta^2q^{\frac{N+1}{2}}$. So
$$|B_{i}^{\prime}(v)|\geq \theta b_iq^{\frac{N+1}{2}}(1 - o(\log^{-10}q)) - 8a_i\theta^2b_i q^{\frac{N+1}{2}} -6a_i\theta^2b_i q^{\frac{N+1}{2}}=$$
$$=\theta b_iq^{\frac{N+1}{2}}(1 - o(\log^{-10}q)) - 14a_i\theta^2b_i q^{\frac{N+1}{2}}.$$

\end{itemize}
Note that since
$$1-P_i^{u}=1-p_i\max_v|A_i(v)|(1+o(1))\geq 1-\frac{\theta}{b_iq^{\frac{N+1}{2}}}\frac{3}{2}a_ib_iq^{\frac{N+1}{2}}\geq 1-\frac{3}{2}\theta a_i,$$
it holds that 
$$\frac{b_i}{b_{i+1}}=\frac{1}{1-P_i^u}\leq \frac{1}{1-\frac{3}{2}\theta a_i}\leq 1+2\theta a_i.$$
First of all observe that
$$|B^{\prime}_i(v)|-|U_i(v)|\leq |M_i^{\prime}(v)|\leq |B_{i}^{\prime}(v)|\leq |B_i(v)|.$$
Then
$$|M_i^{\prime}(v)|\leq |B_i(v)|\leq \theta b_i q^{\frac{N+1}{2}}(1+o(\log^{-10}q))\leq$$
$$\leq \theta b_{i+1} q^{\frac{N+1}{2}}(1+2\theta a_i)(1+o(\log^{-10}q))\leq \theta b_{i+1} q^{\frac{N+1}{2}}(1+3\theta a_i).$$
Also,
$$|M_i^{\prime}(v)|\geq |B^{\prime}_i(v)|-|U_i(v)|\geq\theta b_iq^{\frac{N+1}{2}}(1 - o(\log^{-10}q)) - 14a_i\theta^2b_i q^{\frac{N+1}{2}}-6\theta^2 a_ib_iq^{\frac{N+1}{2}}\geq$$
$$\geq \theta b_iq^{\frac{N+1}{2}}(1  - 15a_i\theta)\geq \theta b_{i+1}q^{\frac{N+1}{2}}(1  - 21a_i\theta).$$
Finally we have to estimate $|L^{\prime}|$. Recall that $L=A_i(v)$, so $A(L,j)=A_i(v,j)$ and $a(L)=\max_{u,v}|A_{i}(u,v)|$. By Remark \ref{ProprValide} we have that, since we are in the second phase,
$$a(L)\leq ib_iq+i\log^{40}q\leq \frac{1}{2}\log^3 q\log^{100} q+\log^3q +\log^{40} q\leq \log^{103}q.$$
Also,
$$\mathbb{E}(|L^{\prime}|)\simeq |A_i(v)|\frac{b_{i+1}}{b_i}\simeq a_ib_{i+1}q^{\frac{N+1}{2}}\geq a_i\log^c q\geq \log^{298} q.$$
Note that for some $Z>0$, $a_i b_{i+1}q^{\frac{N+1}{2}}(1-Z(i-1)\theta^2)\leq \mathbb{E}(|L^{\prime}|)\leq a_i b_{i+1}q^{\frac{N+1}{2}}(1+Z(i-1)\theta^2)$ by induction.
Since $a(L)\leq \log^{103}q \leq\frac{\log^{c-2}q}{2\log^{2(K+5)} q} \leq  \frac{a_ib_iq^{\frac{N+1}{2}}}{2\log^{2(K+5)} q}\leq \frac{|L|}{\log^{2(K+5)} q}$, by Corollary \ref{Corollary 4.8}, with very high probability
$$a_i b_{i+1}q^{\frac{N+1}{2}}(1-Z(i-1)\theta^2)(1-o(\log^{-10}q))\leq |L^{\prime}|\leq a_i b_{i+1}q^{\frac{N+1}{2}}(1+Z(i-1)\theta^2)(1+o(\log^{-10}q)).$$
The lower bound on $|A_{i+1}(v)|$ is then
$$|A_{i+1}(v)|=|L^{\prime}|+|M_i^{\prime}(v)|\geq a_i b_{i+1}q^{\frac{N+1}{2}}(1-Z(i-1)\theta^2)(1-o(\log^{-10}q)+\theta b_{i+1}q^{\frac{N+1}{2}}(1  - 21a_i\theta)\geq$$
$$\geq b_{i+1}q^{\frac{N+1}{2}}(a_i+\theta-Za_ii\theta^2-21a_i\theta^2)\geq b_{i+1}q^{\frac{N+1}{2}}(a_{i+1}-\theta^3-Za_ii\theta^2-21a_i\theta^2),$$
see also \cite[Remark 5.1]{KV2003}; since $a_{i+1}\geq a_i$ and we can consider $Z\geq 22$ we have
$$|A_{i+1}(v)|\geq a_{i+1}b_{i+1}q^{\frac{N+1}{2}}(1-Z(i+1)\theta^2).$$
Also, the upper bound on $|A_{i+1}(v)|$ is
$$|A_{i+1}(v)|=|L^{\prime}|+|M_i^{\prime}(v)|\leq a_i b_{i+1}q^{\frac{N+1}{2}}(1+Z(1-i)\theta^2)(1+o(\log^{-10}q))+\theta b_{i+1}q^{\frac{N+1}{2}}(1  +3a_i\theta)\leq$$
$$b_{i+1}q^{\frac{N+1}{2}}(a_i+\theta+Za_ii\theta^2+3a_i\theta^2)\leq a_{i+1}b_{i+1}q^{\frac{N+1}{2}}(1+Zi\theta^2+7\theta)\leq a_{i+1}b_{i+1}q^{\frac{N+1}{2}}(1+Z(i+1)\theta^2),$$
since, by \cite[Remark 5.1]{KV2003} $a_{i+1}\geq a_i+\theta-(3a_i\theta^2+10\theta^3)\geq a_i+\theta-4a_i\theta^2.$
\endproof



\begin{proposition}\label{Property4bis}
Property \eqref{4bis} holds.
\end{proposition}
\proof
We have to prove that with very high probability
$$(1-\log^{-10} q)^{i+1} q^{N} b_{i+1}\leq |S_{i+1}|\leq (1+\log^{-10} q)^{i+1} q^{N} b_{i+1}.$$

By Remark \ref{ProprValide},
$$|S_{i}(\ell)| \leq  2\log^{c_1} q.$$
Let $L=S_{i}$ and $L^{\prime}=S_{i+1}$, $A(L,j)=\{x \in S_i | (jx) \cap A_i \neq \emptyset\}=A_{i}(j)$. An upper bound for $A(L,j)$ is, by Property \eqref{3bis},
$$a(L)\leq  2a_ib_iq^{\frac{N+1}{2}}\leq 2a_iq^{\frac{N-1}{2}}\log^{c_1} q.$$
By hypotesis $|L|\geq  \frac{1}{2}b_iq^{N}\geq\frac{1}{2}q^{\frac{N-1}{2}}\log^c q$. Then
$$a(L)\leq 2a_iq^{\frac{N-1}{2}}\log^{c_1} q \leq \frac{|L|}{\log^{100}q},$$
and by Lemma \ref{Lemma 4.7} with very high probability
$$\left||L^{\prime}| -\mathbb{E}(|L^{\prime}|) \right| \leq (a(L)|L|)^{1/2} \log^{5} q\leq |L|\log ^{-45}q.$$
Since  $\mathbb{E}(|L^{\prime}|)= |L|(1-P_i^{u})$ (see Remark \ref{RemarkComp}), with very high probability
$$|L|(1-P_i^{u})-|L|\log ^{-45}q \leq |L^{\prime}|\leq  |L|(1-P_i^{u})+|L|\log ^{-45}q.$$
For the first inequality
$$|L|(1-P_i^{u}-\log ^{-45}q)\geq$$
$$\geq (1-\log^{-10} q)^i b_iq^{N}(1-P_i^{u}-\log ^{-45}q)\geq (1-\log^{-10} q)^{i+1} b_{i+1}q^{N}.$$
For the second inequality
$$|L|(1-P_i^{u} +\log^{-45} q)\leq $$
$$\leq(1+\log^{-10}q)^i b_i q^N(1-P_i^{u} +\log^{-45} q)    \leq (1+\log^{-10}q)^{i+1} b_{i+1} q^N.$$
In fact,
$$1-P_i^{u} -\log^{-45} q \geq (1-\log^{-10}q)(1-P_i^{u}) \quad \textrm{and} \quad 1-P_i^{u} +\log^{-45} q \leq (1+\log^{-10}q)(1-P_i^{u})$$
since
$$1-P_i^{u}=\frac{b_{i+1}}{b_{i}}\geq \frac{1}{2}.$$

\endproof

\begin{proposition}\label{Property5bis}
Property \eqref{5bis} holds.
\end{proposition}
\proof
We have to prove that with very high probability
$$|\Omega_{i+1}|\leq (1+\log^{-10} q)^{i+1} q^{N} b^{\prime}_{i+1}.$$
As already observed $b_i^{\prime} \leq 2 b_i$; moreover by Remark \ref{ProprValide},
$$|\Omega_{i}(\ell)| \leq 4\log^{c_1} q.$$
Let $L=\Omega_{i}$ and $L^{\prime}=\Omega_{i+1}$, $A(L,j)=\{x \in \Omega_i | (jx) \cap A_i \neq \emptyset\}$. An upper bound for $A(L,j)$ is $\max_{\ell}|\Omega_{i}(\ell)| \times |A_i|$, then
$$a(L)\leq  8a_iq^{\frac{N-1}{2}}\log^{c_1} q,$$
since $|A_i|\leq 2i \theta q^{\frac{N-1}{2}}$. By hypothesis the uncovered points are at least $|L|\geq  \frac{1}{2}b_iq^{N}\geq\frac{1}{2}q^{\frac{N-1}{2}}\log^c q$ (if not, the algorithm can stop and we obtain a cap of the desired size $O\left(q^{\frac{N-1}{2}}\log^c q\right)$). Then
$$a(L)\leq 8a_iq^{\frac{N-1}{2}}\log^{c_1} q \leq \frac{|L|}{\log^{100}q},$$
and by Lemma \ref{Lemma 4.7} with very high probability
$$\left||L^{\prime}| -\mathbb{E}(|L^{\prime}|) \right| \leq (a(L)|L|)^{1/2} \log^{5} q\leq |L|\log ^{-45}q.$$
Since  $\mathbb{E}(|L^{\prime}|)\leq |L|(1-P_i^{\ell})= |L|\frac{b_{i+1}^{\prime}}{b_i^\prime}$, with very high probability
$$|L^{\prime}|\leq |L|(1-P_i^{\ell})+|L|\log ^{-45}q=|L|(1-P_i^{\ell}+\log ^{-45}q)\leq$$
$$\leq (1+\log^{-10}q)^i q^N b_i^{\prime}(1-P_i^{\ell}+\log ^{-45}q)
\leq   (1+\log^{-10}q)^{i+1} q^Nb_{i+1}^{\prime}.$$
In fact
$$b_i^{\prime}(1-P_i^{\ell}+\log ^{-45}q)
\leq   (1+\log^{-10}q) b_{i+1}^{\prime},$$
since $$1-P_i^{\ell}\geq 1-P_i^{u} \geq\frac{1}{2}.$$
\endproof


\end{document}